\documentclass[a4paper,11pt]{amsart}
\usepackage{graphicx}
\usepackage{mathptmx}
\usepackage{amsmath}
\usepackage{amssymb}
\usepackage{enumitem}
\usepackage{xcolor}
\usepackage{xparse}
\NewDocumentCommand{\eulerian}{omm}
 {%
  \genfrac<>{0pt}{}{#2}{#3}%
  \IfValueT{#1}{_{\!#1}}%
 }

 \newmuskip\pFqmuskip

\newcommand*\pFq[6][8]{%
  \begingroup % only local assignments
  \pFqmuskip=#1mu\relax
  \mathchardef\normalcomma=\mathcode`,
  \mathcode`\,=\string"8000
  \begingroup\lccode`\~=`\,
  \lowercase{\endgroup\let~}\pFqcomma
  {}_{#2}F_{#3}{\left(\genfrac..{0pt}{}{#4}{#5}\bigg|#6\right)}%
  \endgroup
}
\newcommand{\pFqcomma}{{\normalcomma}\mskip\pFqmuskip}

\newtheorem{theorem}{Theorem}
\newtheorem{lemma}[theorem]{Lemma}

\newtheorem{remark}[theorem]{Remark}

\begin{document}

\title[ degenerate power series]{A series transformation formula and related degenerate polynomials}

\author{Taekyun  Kim}
\address{Department of Mathematics, Kwangwoon University, Seoul 139-701, Republic of Korea}
\email{tkkim@kw.ac.kr}

\author{Dae San Kim}
\address{Department of Mathematics, Sogang University, Seoul 121-742, Republic of Korea}
\email{dskim@sogang.ac.kr}

%\author{Hyunseok Lee}
%\address{
%Department of Mathematics, Kwangwoon University, Seoul 139-701, Republic of Korea}
%\email{luciasconstant@kw.ac.kr}


\maketitle

\begin{abstract}
Recently, Boyadzhiev studied a power series whose coefficients are binomial expressions and extended some known formulas involving classical special functions and polynomials. The aim of this paper is to adopt his ideas to express several identities involving `degenerate formal power series' as those including degenerate Stirling numbers of the second kind, degenerate Bell polynomials, degenerate Fubini polynomials and degenerate poly-Bernoulli polynomials.
\end{abstract}

\section{Introduction}

In [3], Boyadzhiev studied a power series whose coefficients are binomial expressions. He extended some known formulas involving classical special functions and polynomials like Hurwitz zeta functions, Euler eta functions, Bernoulli polynomials and Euler polynomials.
He had done this by looking at these formulas from a different perspective, including them in a larger theory and connecting them to the series transformation formulas of Euler [4] and the series transformation formulas considered in [5].
This led naturally to certain asymptotic expansions and gave new proofs of some classical asymptotic expansions. In particular, Boyadzhiev obtained the asymptotic expansion of zeta function. \par

In recent years, studying degenerate versions of some special polynomials and numbers regained interests of some mathematicians, which include the degenerate Bernoulli numbers of the second kind, the degenerate Stirling numbers of both kinds, the degenerate Cauchy numbers, the degenerate Bell numbers and polynomials,
the degenerate complete Bell polynomials and numbers, and so on (see [7,9,10,13,14,16,18] and the references therein). It is remarkable that this study of degenerate versions is not only limited to polynomials and numbers but also extended to transcendental functions like the gamma functions (see [11,12]). They have been studied by various means like combinatorial methods, generating functions, differential equations, umbral calculus, $\lambda$-umbral calculus, $p$-adic analysis, and probability theory. \par

The aim of this paper is to adopt the ideas in [3] and to express several identities involving `degenerate formal power series' as those including degenerate Stirling numbers of the second kind, degenerate Bell polynomials, degenerate Fubini polynomials and degenerate poly-Bernoulli polynomials (see Theorems 2,5,6,9-12). Along the way, we also obtain some related identities which involve the $\lambda$-falling factorials, the degenerate Stirling numbers of both kinds, the degenerate Bernoulli numbers, the degenerate Fubini polynomials and the degenerate Euler polynomials. For the rest of this section, we recall the facts that will be used throughout this paper.

For any $\lambda\in\mathbb{R}$, the degenerate exponential function is defined by
\begin{equation}
e_{\lambda}^{x}(t)=\sum_{n=0}^{\infty}\frac{(x)_{n,\lambda}}{n!}t^{n},\quad (\mathrm{see}\ [7,9,10]),\label{1}
\end{equation}
where the $\lambda$-falling factorials are given by $(x)_{0,\lambda}=1,\ (x)_{n,\lambda}=x(x-\lambda)(x-2\lambda)\cdots(x-(n-1)\lambda),\ (n\ge 1)$. \par
When $x=1$, we write $e_{\lambda}(t)=e_{\lambda}^{1}(t)$. Note that $\displaystyle\lim_{\lambda\rightarrow 0}e_{\lambda}^{x}(t)=e^{xt}\displaystyle$. \par
Let us consider the `degenerate formal power series' which is given by
\begin{equation}
f_{\lambda}(t)=\sum_{k=0}^{\infty}a_{k}(t)_{k,\lambda}=a_{0}+a_{1}(t)_{1,\lambda}+a_{2}(t)_{2,\lambda}+\cdots\in\mathbb{C}[\![t]\!],\label{2}
\end{equation}
where the $a_{i}'s$ are constant complex numbers.
Then we let
\begin{equation}
f(t)=\lim_{\lambda\rightarrow 0}f_{\lambda}(t)=\sum_{k=0}^{\infty}a_{k}t^{n}\in	\mathbb{C}[\![t]\!].\label{3}
\end{equation}
The Carlitz degenerate Bernoulli polynomials are defined by
\begin{equation}
\frac{t}{e_{\lambda}(t)-1}e_{\lambda}^{x}(t)=\sum_{n=0}^{\infty}\beta_{n,\lambda}(x)\frac{t^{n}}{n!}	,\quad(\mathrm{see}\ [6]).\label{4}
\end{equation}
Note that $\displaystyle\lim_{n\rightarrow\infty}\beta_{n,\lambda}=B_{n}(x)$, where $B_{n}(x)$ are Bernoulli polynomials given by
\begin{displaymath}
\frac{t}{e^{t}-1}e^{xt}=\sum_{n=0}^{\infty}B_{n}(x)\frac{t^{n}}{n!},\quad (\mathrm{see}\ [1,14-18]).
\end{displaymath}
Kim-Kim considered the degenerate Stirling numbers of the first kind defined by
\begin{equation}
(x)_{n}=\sum_{k=0}^{n}S_{1,\lambda}(n,k)(x)_{k,\lambda},\quad (n\ge 0),\quad (\mathrm{see}\ [7]),\label{5}
\end{equation}
where $(x)_{0}=1,\ (x)_{n}=x(x-1)\cdots(x-n+1),\ (n\ge 1)$. \par
As the inversion formula of \eqref{5}, the degenerate Stirling numbers of the second kind are given by
\begin{equation}
(x)_{n,\lambda}=\sum_{k=0}^{n}S_{2,\lambda}(n,k)(x)_{k},\quad (n\ge 0),\quad (\mathrm{see}\ [7,9,10]).\label{6}
\end{equation}
In [13], the degenerate Bell polynomials are defined by
\begin{align}
e^{x(e_{\lambda}(t)-1)}=\sum_{n=0}^{\infty}\phi_{n,\lambda}(x)\frac{t^{n}}{n!}. \label{7}
\end{align}
Note that $\displaystyle \phi_{n,\lambda}(x)=\sum_{k=0}^{n}S_{2,\lambda}(n,k)x^{k}\displaystyle$. When $x=1,\ \phi_{n,\lambda}=\phi_{n,\lambda}(1)$ are called the degenerate Bell numbers. \par
For any $\lambda\in\mathbb{R}$, the degenerate Fubini polynomials are defined by Kim-Kim as
\begin{equation}
\frac{1}{1-x(e_{\lambda}(t)-1)}=\sum_{n=0}^{\infty}F_{n,\lambda}(x)\frac{t^{n}}{n!},\quad (\mathrm{see}\ [14]).\label{8}	
\end{equation}
It is well known that the polylogarithmic function is defined by
\begin{equation}
\mathrm{Li}_{k}(x)=\sum_{n=1}^{\infty}\frac{x^{k}}{n^{k}},\quad (k\in\mathbb{Z}),\quad (\mathrm{see}\ [2]).\label{9}
\end{equation}
Note that $\mathrm{Li}_{1}(x)=-\log (1-x)$.\par

\section{Some identities of degenerate special polynomials}
Here we consider the `degenerate formal power series' given by
\begin{align}
f_{\lambda}(t)&=\sum_{k=0}^{\infty}a_{k}(t)_{k,\lambda}\label{10}	\\
&=a_{0}+a_{1}(t)_{1,\lambda}+a_{2}(t)_{2,\lambda}+a_{3}(t)_{3,\lambda}+\cdots\in\mathbb{C}[\![t]\!]\nonumber,
\end{align}
and express several identities involving `degenerate formal power series' as those including degenerate Stirling numbers of the second kind, degenerate Bell polynomials, degenerate Fubini polynomials and degenerate poly-Bernoulli polynomials (see Theorems 2,5,6,9-12).
From \eqref{6}, we note that
\begin{align}
\sum_{n=k}^{\infty}S_{2,\lambda}(n,k)\frac{t^{n}}{n!}&=\frac{1}{k!}\big(e_{\lambda}(t)-1\big)^{k}=\frac{1}{k!}\sum_{l=0}^{k}\binom{k}{l}(-1)^{k-l}e_{\lambda}^{l}(t)\label{11}\\
&=\frac{1}{k!}\sum_{l=0}^{k}\binom{k}{l}(-1)^{k-l}\sum_{n=0}^{\infty}(l)_{n,\lambda}\frac{t^{n}}{n!}\nonumber \\
&=\sum_{n=0}^{\infty}\bigg\{\frac{1}{k!}\sum_{l=0}^{k}\binom{k}{l}(-1)^{k-l}(l)_{n,\lambda}\bigg\}
\frac{t^{n}}{n!}.\nonumber
\end{align}
By comparing the coefficients on both sides of \eqref{11}, we get
\begin{equation}
\frac{1}{k!}\sum_{l=0}^{k}\binom{k}{l}(-1)^{k-l}(l)_{n,\lambda}=S_{2,\lambda}(n,k), \label{12}
\end{equation}
with the understanding that $S_{2,\lambda}(n,k)=0$,~ for $0 \le n < k$.
It is not difficult to show that
\begin{equation}
(y+zk)_{m,\lambda}=\sum_{l=0}^{m}\binom{m}{l}(zk)_{l,\lambda}(y)_{m-l,\lambda},\quad (m\ge 0).\label{13}
\end{equation}
We observe that
\begin{align}
(zk)_{p,\lambda}&=(zk)(zk-\lambda)\cdots(zk-(p-1)\lambda)\label{14}\\
&=k^{p}z\bigg(z-\frac{\lambda}{k}\bigg)\cdots\bigg(z-(p-1)\frac{\lambda}{k}\bigg)	\nonumber \\
&=k^{p}(z)_{p,\frac{\lambda}{k}}=z^{p} (k)_{p,\frac{\lambda}{z}},\quad (p\ge 0).\nonumber
\end{align}
From \eqref{13} and \eqref{14}, we note that
\begin{align}
(y+zk)_{m,\lambda}&=\sum_{p=0}^{m}\binom{m}{p}k^{p}(z)_{p,\frac{\lambda}{k}}(y)_{m-p,\lambda}
\label{15}\\
&=\sum_{p=0}^{m}\binom{m}{p}z^{p}	(y)_{m-p,\lambda} (k)_{p,\frac{\lambda}{z}}.\nonumber
\end{align}
By \eqref{12} and \eqref{15}, we get
\begin{align}
\sum_{k=0}^{n}\binom{n}{k}(-1)^{k} (y+zk)_{m,\lambda}&=\sum_{k=0}^{n}\binom{n}{k}(-1)^{k}\bigg\{\sum_{p=0}^{m}\binom{m}{p}z^{p}(y)_{m-p,\lambda} (k)_{p,\frac{\lambda}{z}} \bigg\}\label{16} \\
&=\sum_{p=0}^{m}\binom{m}{p}z^{p}(y)_{m-p,\lambda}\bigg\{\sum_{k=0}^{n}\binom{n}{k}(-1)^{k}(k)_{p,\frac{\lambda}{z}}\bigg\}\nonumber \\
&=(-1)^{n}n!\sum_{p=0}^{m}\binom{m}{p}z^{p}(y)_{m-p,\lambda}\bigg(\frac{1}{n!}\sum_{k=0}^{n}\binom{n}{k}(-1)^{n-k} (k)_{p,\frac{\lambda}{z}}\bigg) \nonumber \\
&=(-1)^{n}n!\sum_{p=0}^{m}\binom{m}{p}z^{p}(y)_{m-p,\lambda}S_{2,\frac{\lambda}{z}}(p,n).\nonumber
\end{align}
Therefore, by comparing the coefficients on both sides of \eqref{16}, we obtain the following theorem.
\begin{theorem}
For $n\ge 0$, the following identity holds.
\begin{displaymath}
\sum_{k=0}^{n}\binom{n}{k}(-1)^{k} (y+zk)_{m,\lambda}= (-1)^{n}n!\sum_{p=0}^{m}\binom{m}{p}z^{p}(y)_{m-p,\lambda}S_{2,\frac{\lambda}{z}}(p,n).
\end{displaymath}
\end{theorem}
\begin{remark}
As $S_{2,\lambda}(n,k)=0$, for $k > n$ and $k < 0 $, we see from Theorem 1 that we have
\begin{align*}
\sum_{k=0}^{n}\binom{n}{k}(-1)^{k} (y+zk)_{m,\lambda}=0, ~~\mathrm{if}~~ n > m.
\end{align*}
\end{remark}
From \eqref{10}, we note that
\begin{equation}
f_{\lambda}(y+zk)=\sum_{m=0}^{\infty}a_{m}(y+zk)_{m,\lambda}.\label{17}	
\end{equation}
By \eqref{12}, \eqref{15}, \eqref{17} and Theorem 1, we get
\begin{align}
\sum_{k=0}^{n}\binom{n}{k}(-1)^{k}f_{\lambda}(y+zk)&= \sum_{k=0}^{n}\binom{n}{k}(-1)^{k}\sum_{m=0}^{\infty}a_{m}(y+zk)_{m,\lambda}\label{18}	\\
&=\sum_{m=0}^{\infty}a_{m}(-1)^{n}n!\sum_{p=0}^{m}\binom{m}{p}S_{2,\frac{\lambda}{z}}(p,m)z^{p}(y)_{m-p,\lambda}\nonumber\\
&=(-1)^{n}n!\sum_{m=0}^{\infty}a_{m}\bigg\{\sum_{p=0}^{m}\binom{m}{p}S_{2,\frac{\lambda}{z}}(p,n)z^{p}(y)_{m-p,\lambda}\bigg\}.\nonumber
\end{align}
Therefore, we obtain the following theorem.
\begin{theorem}
For $n\ge 0$, the following identity is valid.
\begin{displaymath}
\sum_{k=0}^{n}\binom{n}{k}(-1)^{k}f_{\lambda}(y+zk)=(-1)^{n}n!\sum_{m=0}^{\infty}a_{m}\bigg\{\sum_{p=0}^{m}\binom{m}{p}S_{2,\frac{\lambda}{z}}(p,n)z^{p}(y)_{m-p,\lambda}\bigg\}.
\end{displaymath}
In particular, for $y=0$,
\begin{displaymath}
\sum_{k=0}^{n}\binom{n}{k}(-1)^{k}f_{\lambda}(zk)=(-1)^{n}n!\sum_{m=0}^{\infty}a_{m}S_{2,\frac{\lambda}{z}}(m,n)z^{m}.
\end{displaymath}
\end{theorem}
If $n=m$ in Theorem 1, then we have
\begin{align}
\sum_{k=0}^{n}\binom{n}{k}(-1)^{k}(y+zk)_{n,\lambda}&=(-1)^{n}n!\sum_{p=0}^{n}\binom{n}{p}z^{p}(y)_{n-p,\lambda}S_{2,\frac{\lambda}{z}}(p,n)\label{20} \\
&=(-1)^{n}n!z^{n}.\nonumber
\end{align}
Therefore, by \eqref{20}, we obtain the following theorem.
\begin{theorem}
For $n\ge 0$, we have the following identity:
\begin{displaymath}
\sum_{k=0}^{n}\binom{n}{k}(-1)^{k}(y+zk)_{n,\lambda}= (-1)^{n}n!z^{n}.
\end{displaymath}
\end{theorem}
For $n\in\mathbb{N}$, we have
\begin{align}
(x)_{n,\lambda}=\sum_{k=0}^{n}S_{2,\lambda}(n,k)(x)_{k}=\sum_{k=1}^{n}S_{2,\lambda}(n,k)(x)_{k}.\label{21}
\end{align}
Dividing \eqref{21} by $x$ and then taking $x=0$, we have
\begin{equation}
(-1)^{n-1}\lambda^{n-1}(n-1)!=\sum_{k=1}^{n}S_{2,\lambda}(n,k)(-1)^{k-1}(k-1)!.\label{22}
\end{equation}
Therefore, by \eqref{22}, we obtain the following lemma.
\begin{lemma}
For $n\in\mathbb{N}$, the following identity holds true.
\begin{displaymath}
(-1)^{n-1}\lambda^{n-1}(n-1)!=\sum_{k=1}^{n}S_{2,\lambda}(n,k)(-1)^{k-1}(k-1)!.
\end{displaymath}	
\end{lemma}
From Theorem 3, we note that
\begin{displaymath}
\sum_{k=0}^{n}\binom{n}{k}(-1)^{k}f_{\lambda}(zk)=(-1)^{n}n!\sum_{m=0}^{\infty}a_{m}z^{m}S_{2,\frac{\lambda}{z}}(m,n).
\end{displaymath}
By multiplying $\frac{1}{n}$ on both sides, we get
\begin{displaymath}
\frac{1}{n}\sum_{k=0}^{n}\binom{n}{k}(-1)^{k}f_{\lambda}(zk)=(-1)^{n}(n-1)!\sum_{m=0}^{\infty}a_{m}z^{m}S_{2,\frac{\lambda}{z}}(m,n),
\end{displaymath}
where $n=1,2,3,\dots$. \par
Thus, by using Lemma 5, we have
\begin{align}
\sum_{n=1}^{\infty}\frac{1}{n}\bigg\{\sum_{k=0}^{n}\binom{n}{k}(-1)^{k}f_{\lambda}(zk)\bigg\}&=\sum_{n=1}^{\infty}(-1)^{n}(n-1)!	\sum_{m=0}^{\infty}a_{m}z^{m}S_{2,\frac{\lambda}{z}}(m,n)\label{23} \\
&=\sum_{m=1}^{\infty}a_{m}z^{m}\sum_{n=1}^{\infty}(-1)^{n}(n-1)!S_{2,\frac{\lambda}{z}}(m,n)\nonumber \\
&=\sum_{m=1}^{\infty}a_{m}z^{m}(-1)^{m}\bigg(\frac{\lambda}{z}\bigg)^{m-1}(m-1)! \nonumber \\
&=z\sum_{m=1}^{\infty}a_{m}(-1)^{m}\lambda^{m-1}(m-1)!.\nonumber
\end{align}
Therefore, by \eqref{23}, we obtain the following theorem.
\begin{theorem}
For $\displaystyle f_{\lambda}(t)=\sum_{k=0}^{\infty}a_{k}(t)_{k,\lambda}\in\mathbb{C}[\![t]\!]\displaystyle$, the following identity holds.
\begin{displaymath}
\sum_{n=1}^{\infty}\frac{1}{n}\bigg\{\sum_{k=0}^{n}\binom{n}{k}(-1)^{k}f_{\lambda}(zk)\bigg\}=-z\sum_{m=1}^{\infty}a_{m}(-1)^{m}\lambda^{m-1}(m-1)!
\end{displaymath}
\end{theorem}
We observe that
\begin{align}
f_{\lambda}^{\prime}(0)=\frac{d}{dt}f_{\lambda}(t)\bigg|_{t=0}=\sum_{k=1}^{\infty}a_{k}\frac{d}{dt}(t)_{k,\lambda}\bigg|_{t=0}=\sum_{k=1}^{\infty}a_{k}(-1)^{k-1}\lambda^{k-1}(k-1)!.\label{24}
\end{align}
By Theorem 6 and \eqref{24}, we get
\begin{align}
\sum_{n=1}^{\infty}\frac{1}{n}\bigg\{\sum_{k=0}^{n}\binom{n}{k}(-1)^{k}f_{\lambda}(zk)\bigg\}=-z\sum_{m=1}^{\infty}a_{m}(-1)^{m-1}\lambda^{m-1}(m-1)!=-zf_{\lambda}^{\prime}(0).\label{25}	
\end{align}
Therefore, we obtain the following theorem.
\begin{theorem}
Let $\displaystyle f_{\lambda}(t)=\sum_{k=0}^{\infty}a_{k}(t)_{k,\lambda}\in\mathbb{C}[\![t]\!]\displaystyle$. Then the following identity is valid.
\begin{displaymath}
\sum_{n=1}^{\infty}\frac{1}{n}\bigg\{\sum_{k=0}^{n}\binom{n}{k}(-1)^{k}f_{\lambda}(zk)\bigg\}=-zf_{\lambda}^{\prime}(0),
\end{displaymath}
where $f_{\lambda}^{\prime}(0)=\frac{d}{dt}f_{\lambda}(t)\big|_{t=0}$.
\end{theorem}
Let us take $f_{\lambda}(t)=\binom{t}{p}_{\lambda}=\frac{(t)_{p,\lambda}}{p!},\ (p\in\mathbb{N}\cup\{0\})$. Then we have $a_{m}=\frac{1}{m!}\delta_{m,p},\ (m,p\ge 0)$, where $\delta_{m,n}$ is the Kronecker's symbol. \\
Then, from Theorem 3, we have
\begin{align}
\sum_{k=0}^{n}\binom{n}{k}(-1)^{k}\binom{zk}{p}_{\lambda}
=(-1)^{n}\frac{n!}{p!}S_{2,\frac{\lambda}{z}}(p,n)z^{p}.\label{26}
\end{align}
Thus, we have the following lemma.
\begin{lemma}
For $n\ge 0$, we have the following identity.
\begin{equation}
\sum_{k=0}^{n}\binom{n}{k}(-1)^{k}\binom{zk}{p}_{\lambda}= (-1)^{n}\frac{n!}{p!}S_{2,\frac{\lambda}{z}}(p,n)z^{p}. \label{27}
\end{equation}
\end{lemma}
Note that $\displaystyle\sum_{k=0}^{n}\binom{n}{k}(-1)^{k}\binom{zk}{n}_{\lambda}=(-1)^{n}z^{n}\displaystyle$. It is well known that the Bell polynomials are defined by
\begin{displaymath}
\sum_{n=0}^{\infty}\phi_{n}(x)\frac{t^{n}}{n!}=\lim_{\lambda\rightarrow 0}e^{x(e_{\lambda}(t)-1)}=e^{x(e^{t}-1)}.
\end{displaymath}
The compositional inverse of $e_{\lambda}(t)$ is denoted by $\log_{\lambda}t$. So $e_{\lambda}(\log_{\lambda} t)=\log_{\lambda}(e_{\lambda}(t))=t$. They are called the degenerate logarithms and given by
\begin{equation}
\log_{\lambda}t=\frac{1}{\lambda}(t^{\lambda}-1),~~\log_{\lambda}(1+t)=\frac{1}{\lambda}\big((1+t)^{\lambda}-1\big)=\sum_{n=1}^{\infty}\lambda^{n-1}(1)_{n,\frac{1}{\lambda}}\frac{t^{n}}{n!}.\label{28}
\end{equation}
Now, we observe that
\begin{align}
\frac{1}{n!}\Big((1+t)^{z}-1\Big)^{n}&=\frac{1}{n!}\Big(e^{z\log (1+t)}-1\Big)^{n}\label{29}\\
&=\sum_{m=n}^{\infty}S_{2}(m,n)\frac{z^{m}}{m!}\big(\log(1+t)\big)^{m}=\sum_{m=n}^{\infty}S_{2}(m,n)z^{m}\sum_{p=m}^{\infty}S_{1}(p,m)\frac{t^{p}}{p!}\nonumber \\
&=\sum_{p=n}^{\infty}\bigg(\sum_{m=n}^{p}S_{2}(m,n)S_{1}(p,m)z^{m}\bigg)\frac{t^{n}}{p!},\nonumber
\end{align}
where $S_{2}(m,n)=\lim_{\lambda\rightarrow 0}S_{2,\lambda}(m,n)$ and $S_{1}(m,n)=\lim_{\lambda\rightarrow 0}S_{1,\lambda}(m,n)$ are respectively the Stirling numbers of the second kind and the Stirling numbers of the first kind. \par
On the other hand,
\begin{align}
\frac{1}{n!}\Big((1+t)^{z}-1\Big)^{n}&=\frac{1}{n!}\Big(e_{\lambda}^{z}(\log_{\lambda}(1+t))-1\Big)^{n}\label{30}\\
&=\frac{1}{n!}\Big(e_{\frac{\lambda}{z}}(z\log_{\lambda}(1+t))-1\Big)^{n}\nonumber\\
&=\sum_{m=n}^{\infty}S_{2,\frac{\lambda}{z}}(m,n)\frac{z^{m}}{m!}\big(\log_{\lambda}(1+t)\big)^{m}\nonumber\\
&=\sum_{m=n}^{\infty}S_{2,\frac{\lambda}{z}}(m,n)z^{m}\sum_{p=m}^{\infty}S_{1,\lambda}(p,m)\frac{t^{p}}{p!}\nonumber \\
&=\sum_{p=n}^{\infty}\bigg(\sum_{m=n}^{p}S_{2,\frac{\lambda}{z}}(m,n)S_{1,\lambda}(p,m)z^{m}\bigg)\frac{t^{p}}{p!}.\nonumber
\end{align}
Therefore, by \eqref{28} and \eqref{30}, we obtain the following theorem.
\begin{theorem}
For $n,p$, with $n \ge p \ge 0$, the following holds true.
\begin{displaymath}
\sum_{m=n}^{p}S_{2}(m,n)S_{1}(p,m)z^{m}=\sum_{m=n}^{p}S_{2,\frac{\lambda}{z}}(m,n)S_{1,\lambda}(p,m)z^{m}.
\end{displaymath}	
\end{theorem}
From Theorem 3, we note that
\begin{align}
&\sum_{n=0}^{\infty}\frac{x^{n}}{n!}\bigg\{\sum_{k=0}^{n}\binom{n}{k}(-1)^{k}f_{\lambda}(y+kz)	\bigg\}\ \label{31}\\
&=\sum_{n=0}^{\infty}\frac{x^{n}}{n!}(-1)^{n}n!\sum_{m=0}^{\infty}a_{m}\bigg\{\sum_{p=0}^{m}\binom{m}{p}S_{2,\frac{\lambda}{z}}(p,n)z^{p}(y)_{m-p,\lambda}\bigg\}\nonumber \\
&=\sum_{m=0}^{\infty}a_{m}\sum_{p=0}^{m}\binom{m}{p}(y)_{m-p,\lambda}z^{p}\bigg\{\sum_{n=0}^{\infty}x^{n}(-1)^{n}S_{2,\frac{\lambda}{z}}(p,n)\bigg\}\nonumber \\
&=\sum_{m=0}^{\infty}a_{m}\sum_{p=0}^{m}\binom{m}{p}(y)_{m-p,\lambda}z^{p}\sum_{n=0}^{p}x^{n}(-1)^{n}S_{2,\frac{\lambda}{z}}(p,n) \nonumber \\
&=\sum_{m=0}^{\infty}a_{m}\sum_{p=0}^{m}\binom{m}{p}(y)_{m-p,\lambda}z^{p}\phi_{p,\frac{\lambda}{z}}(-x).\nonumber
\end{align}
From \eqref{31}, we obtain the following theorem.
\begin{theorem}
Let $\displaystyle f_{\lambda}(t)=\sum_{k=0}^{\infty}a_{k}(t)_{k,\lambda}\in\mathbb{C}[\![t]\!]\displaystyle$. Then we have the following identity:
\begin{equation*}
\sum_{n=0}^{\infty}\frac{x^{n}}{n!}\bigg\{\sum_{k=0}^{n}\binom{n}{k}(-1)^{k}f_{\lambda}(y+kz)\bigg\}= \sum_{m=0}^{\infty}a_{m}\sum_{p=0}^{m}\binom{m}{p}(y)_{m-p,\lambda}z^{p}\phi_{p,\frac{\lambda}{z}}(-x).
\end{equation*}	
In particular, for $y=0$,
\begin{equation}
\sum_{n=0}^{\infty}\frac{x^{n}}{n!}\bigg\{\sum_{k=0}^{n}\binom{n}{k}(-1)^{k}f_{\lambda}(kz)	\bigg\}= \sum_{m=0}^{\infty}a_{m}z^{m}\phi_{m,\frac{\lambda}{z}}(-x).\label{32}
\end{equation}	
\end{theorem}
Let us take $x=-1$ in \eqref{32}. Then we have
\begin{equation}
\sum_{n=0}^{\infty}\frac{(-1)^{n}}{n!}\bigg\{\sum_{k=0}^{n}\binom{n}{k}(-1)^{k}f_{\lambda}(kz)	\bigg\}= \sum_{m=0}^{\infty}a_{m}z^{m}\phi_{m,\frac{\lambda}{z}}.	\label{33}
\end{equation}
Let us take $f_{\lambda}(t)=\binom{t}{p}_{\lambda}=\frac{(t)_{p,\lambda}}{p!}$. Then, from \eqref{32} and Theorem 3, we get
\begin{align}
&\sum_{n=0}^{p}\frac{x^{n}}{n!}\bigg\{\sum_{k=0}^{n}\binom{n}{k}(-1)^{k}\binom{zk}{p}_{\lambda}\bigg\}=\frac{z^{p}}{p!}\phi_{p,\frac{\lambda}{z}}(-x).\label{34}
\end{align}
Here we observe that
$\sum_{k=0}^{n}\binom{n}{k}(-1)^{k}\binom{zk}{p}_{\lambda}=\frac{(-1)^{n}n!}{p!}S_{2,\frac{\lambda}{z}}(p,n)z^{p}=0$,~~if $n > p$.

From \eqref{8}, we note that
\begin{align}
\sum_{n=0}^{\infty}F_{n,\lambda}(x)\frac{t^{n}}{n!}&=\frac{1}{1-x(e_{\lambda}(t)-1)}=\sum_{k=0}^{\infty}x^{k}(e_{\lambda}(t)-1)^{k}\label{35}\\
&=\sum_{k=0}^{\infty}x^{k}k!\frac{1}{k!}(e_{\lambda}(t)-1)^{k}=\sum_{k=0}^{\infty}x^{k}k!	\sum_{n=k}^{\infty}S_{2,\lambda}(n,k)\frac{t^{n}}{n!}\nonumber \\
&=\sum_{n=0}^{\infty}\bigg(\sum_{k=0}^{n}x^{k}k!S_{2,\lambda}(n,k)\bigg)\frac{t^{n}}{n!}.\nonumber
\end{align}
Thus, we have
\begin{equation}
F_{n,\lambda}(x)=\sum_{k=0}^{n}k!S_{2,\lambda}(n,k)x^{k},\quad (n\ge 0).\label{36}	
\end{equation}
From Theorem 3, we have
\begin{align}
&\sum_{n=0}^{\infty}x^{n}\bigg\{\sum_{k=0}^{n}\binom{n}{k}(-1)^{k}f_{\lambda}(y+zk)\bigg\}\label{37}\\
&=\sum_{n=0}^{\infty}x^{n}(-1)^{n}n!\sum_{m=0}^{\infty}a_{m}\bigg\{\sum_{p=0}^{m}\binom{m}{p}(y)_{m-p,\lambda}z^{p}S_{2,\frac{\lambda}{z}}(p,n)\bigg\}\nonumber \\
&=\sum_{m=0}^{\infty}a_{m}\sum_{p=0}^{m}\binom{m}{p}(y)_{m-p,\lambda}z^{p}\bigg\{\sum_{n=0}^{\infty}x^{n}(-1)^{n}n!
S_{2,\frac{\lambda}{z}}(p,n)\bigg\}\nonumber \\
&=\sum_{m=0}^{\infty}a_{m}\sum_{p=0}^{m}\binom{m}{p}(y)_{m-p,\lambda}z^{p}\bigg\{\sum_{n=0}^{p}x^{n}(-1)^{n}n!S_{2,\frac{\lambda}{z}}(p,n)\bigg\}\nonumber \\
&=\sum_{m=0}^{\infty}a_{m}\sum_{p=0}^{m}\binom{m}{p}(y)_{m-p,\lambda}z^{p}F_{p,\frac{\lambda}{z}}(-x).\nonumber
\end{align}
Therefore, by \eqref{37}, we obtain the following theorem.
\begin{theorem}
Let 	$\displaystyle f_{\lambda}(t)=\sum_{k=0}^{\infty}a_{k}(t)_{k,\lambda}\in\mathbb{C}[\![t]\!]\displaystyle$. Then the following identity holds true.
\begin{displaymath}
\sum_{n=0}^{\infty}x^{n}\bigg\{\sum_{k=0}^{n}\binom{n}{k}(-1)^{k}f_{\lambda}(y+zk)\bigg\}= \sum_{m=0}^{\infty}a_{m}\sum_{p=0}^{m}\binom{m}{p}(y)_{m-p,\lambda}z^{p}F_{p,\frac{\lambda}{z}}(-x).
\end{displaymath}
In particular, for $y=0$,
\begin{displaymath}
\sum_{n=0}^{\infty}x^{n}\bigg\{\sum_{k=0}^{n}\binom{n}{k}(-1)^{k}f_{\lambda}(zk)\bigg\}= \sum_{m=0}^{\infty}a_{m}z^{m}F_{m,\frac{\lambda}{z}}(-x).
\end{displaymath}
\end{theorem}
Let us take $f_{\lambda}(t)=\binom{t}{p}_{\lambda}=\frac{(t)_{p,\lambda}}{p!}$. Then we have $a_{m}=\frac{1}{m!}\delta_{m,p},\ (m,p\ge 0)$. Hence we have
\begin{displaymath}
\sum_{n=0}^{\infty}x^{n}\bigg\{\sum_{k=0}^{n}\binom{n}{k}(-1)^{k}\binom{zk}{p}_{\lambda}\bigg\}=\frac{1}{p!}z^{p}F_{p,\frac{\lambda}{z}}(-x).
\end{displaymath}
From Theorem 11, we note that
\begin{equation}
\int_{0}^{x}\sum_{n=0}^{\infty}t^{n}\bigg\{\sum_{k=0}^{n}\binom{n}{k}(-1)^{k}f_{\lambda}(y+zk)\bigg\}dt=\sum_{m=0}^{\infty}a_{m}\sum_{p=0}^{m}\binom{m}{p}(y)_{m-p,\lambda}z^{p}\int_{0}^{x}F_{p,\frac{\lambda}{z}}(-t)dt.\label{38}
\end{equation}
By \eqref{36}, we easily get
\begin{equation}
\int_{0}^{x}F_{p,\frac{\lambda}{z}}(-t)dt=\sum_{j=0}^{p}(-1)^{j}j!\frac{x^{j+1}}{j+1}S_{2,\frac{\lambda}{z}}(p,j).\label{39}
\end{equation}
From \eqref{38} and \eqref{39}, we have
\begin{align}
&\sum_{n=0}^{\infty}\frac{x^{n+1}}{n+1}\bigg\{\sum_{k=0}^{n}\binom{n}{k}(-1)^{k}f_{\lambda}(y+zk)\bigg\} \label{40}\\
&=\sum_{m=0}^{\infty}a_{m}\sum_{p=0}^{m}\binom{m}{p}(y)_{m-p,\lambda}z^{p}\sum_{j=0}^{p}(-1)^{j}j!\frac{x^{j+1}}{j+1}S_{2,\frac{\lambda}{z}}(p,j).\nonumber
\end{align}
Thus, we see that
\begin{align}
&\sum_{n=0}^{\infty}\frac{x^{n}}{n+1}\bigg\{\sum_{k=0}^{n}\binom{n}{k}(-1)^{k}f_{\lambda}(y+zk)\bigg\}\label{41}\\
&=\sum_{m=0}^{\infty}a_{m}\sum_{p=0}^{m}\binom{m}{p}(y)_{m-p,\lambda}z^{p}\sum_{j=0}^{p}(-1)^{j}j!\frac{x^{j}}{j+1}S_{2,\frac{\lambda}{z}}(p,j).\nonumber
\end{align}
Now, we consider $(r-1)$-times iterated integration with respect to $x$ as follows:
\begin{align}
&\frac{1}{x}\underbrace{\int_{0}^{x}\frac{1}{x}\int_{0}^{x}\cdots\frac{1}{x}\int_{0}^{x}}_{(r-1)-\mathrm{times}}\sum_{n=0}^{\infty}\frac{x^{n}}{n+1}\bigg\{\sum_{k=0}^{n}\binom{n}{k}(-1)^{k}f_{\lambda}(y+zk)\bigg\}dx\cdots dx\label{42}\\
&=\sum_{m=0}a_{m}\sum_{p=0}^{m}\binom{m}{p}(y)_{m-p,\lambda}z^{p}\sum_{j=0}^{p}(-1)^{j}\frac{j!}{j+1}S_{2,\frac{\lambda}{z}}(p,j)\nonumber\\
&\quad\quad\quad\times \frac{1}{x}\underbrace{\int_{0}^{x}\frac{1}{x}\int_{0}^{x}\cdots\frac{1}{x}\int_{0}^{x}}_{(r-1)-\mathrm{times}}x^{j}dx\cdots dx\nonumber \\
&=\sum_{m=0}^{\infty}a_{m}\sum_{p=0}^{m}\binom{m}{p}(y)_{m-p,\lambda}z^{p}\sum_{j=0}^{p}(-1)^{j}	\frac{j!}{(j+1)^{r}}x^{j}S_{2,\frac{\lambda}{z}}(p,j).\nonumber
\end{align}
By \eqref{42}, we get
\begin{align}
&\sum_{n=0}^{\infty}\frac{x^{n}}{(n+1)^{r}}\bigg\{\sum_{k=0}^{n}\binom{n}{k}(-1)^{k}f_{\lambda}(y+zk)\bigg\}	\label{43} \\
&=\sum_{m=0}^{\infty}a_{m}\bigg\{\sum_{p=0}^{m}\binom{m}{p}(y)_{m-p,\lambda}z^{p}\sum_{j=0}^{p}(-1)^{j}j!
S_{2,\frac{\lambda}{z}}(p,j)\bigg\}\frac{x^{j}}{(j+1)^{r}}\nonumber
\end{align}
Let us take $y=0$ in \eqref{43}. Then we have
\begin{equation}
\sum_{n=0}^{\infty}\frac{x^{n}}{(n+1)^{r}}\bigg\{\sum_{k=0}^{n}\binom{n}{k}(-1)^{k}f_{\lambda}(zk)\bigg\}=\sum_{m=0}^{\infty}a_{m}z^{m}\sum_{j=0}^{m}(-1)^{j}j!S_{2,\frac{\lambda}{z}}(m,j)\frac{x^{j}}{(j+1)^{r}}.\label{44}	
\end{equation}
By Theorem 11 and \eqref{36}, we see that
\begin{align}
&\sum_{n=0}^{\infty}x^{n}\bigg\{\sum_{k=0}^{n}\binom{n}{k}(-1)^{k}f_{\lambda}(zk)\bigg\}	\label{45} \\
&=\sum_{m=0}^{\infty}a_{m}z^{m}F_{m,\frac{\lambda}{z}}(-x)=\sum_{m=0}^{\infty}a_{m}z^{m}\bigg\{\sum_{p=0}^{m}S_{2,\frac{\lambda}{z}}(m,p)(-x)^{p}p!\bigg\}\nonumber .
\end{align}
Thus, we have
\begin{equation}
\sum_{n=1}^{\infty}x^{n-1}\bigg\{\sum_{k=0}^{n}\binom{n}{k}(-1)^{k}f_{\lambda}(zk)\bigg\}=\sum_{m=1}^{\infty}a_{m}z^{m}\sum_{p=1}^{m}S_{2,\frac{\lambda}{z}}(m,p)(-1)^{p}p!x^{p-1}.\label{46}	
\end{equation}
From \eqref{46}, we note that
\begin{align}
\sum_{n=1}^{\infty}\frac{x^{n}}{n}\bigg\{\sum_{k=0}^{n}\binom{n}{k}(-1)^{k}f_{\lambda}(zk)\bigg\}&=\sum_{m=1}^{\infty}a_{m}z^{m}\sum_{p=1}^{m}S_{2,\frac{\lambda}{z}}(m,p)(-1)^{p}p!\int_{0}^{x}x^{p-1}dx  \label{47}\\
&=\sum_{m=1}^{\infty}a_{m}z^{m}\sum_{p=1}^{m}S_{2,\frac{\lambda}{z}}(m,p)(-1)^{p}\frac{p!}{p}x^{p}.\nonumber
\end{align}
Therefore, we obtain the following theorem.
\begin{theorem}
Let $\displaystyle f_{\lambda}(t)=\sum_{k=0}^{\infty}a_{k}(t)_{k,\lambda}\in\mathbb{C}[\![t]\!]\displaystyle$. Then the following identity holds.
\begin{displaymath}
\sum_{n=1}^{\infty}\frac{x^{n}}{n}\bigg\{\sum_{k=0}^{n}\binom{n}{k}(-1)^{k}f_{\lambda}(zk)\bigg\}=\sum_{m=1}^{\infty}a_{m}z^{m}\sum_{p=1}^{m}S_{2,\frac{\lambda}{z}}(m,p)(-1)^{p}(p-1)!x^{p}.
\end{displaymath}	
\end{theorem}
We define the degenerate poly-Bernoulli polynomials which are given by
\begin{equation}
\frac{\mathrm{Li}_{k}(1-e_{\lambda}(-t))}{1-e_{\lambda}(-t)}e_{\lambda}^{x}(-t)=\sum_{n=0}^{\infty}\beta_{n,\lambda}^{(k)}(x)\frac{t^{n}}{n!}.\label{48}
\end{equation}
When $x=0$, $\beta_{n,\lambda}^{(k)}=\beta_{n,\lambda}^{(k)}(0)$ are called the degenerate poly-Bernoulli numbers. \par
By \eqref{48}, we see that
\begin{align}
\sum_{n=0}^{\infty}\beta_{n,\lambda}^{(k)}\frac{t^{n}}{n!}&=\frac{1}{1-e_{\lambda}(-t)}\sum_{j=1}^{\infty}\frac{(1-e_{\lambda}(-t))^{j}}{j^{k}}=\sum_{j=1}^{\infty}\frac{(1-e_{\lambda}(-t))^{j-1}}{j^{k}}\label{49} \\
&=\sum_{j=0}^{\infty}\frac{(-1)^{j}j!}{(j+1)^{k}}\frac{1}{j!}\big(e_{\lambda}(-t)-1\big)^{j}=\sum_{n=0}^{\infty}\bigg((-1)^{n}\sum_{j=0}^{n}\frac{(-1)^{j}j!}{(j+1)^{k}}S_{2,\lambda}(n,j)\bigg)\frac{t^{n}}{n!}.\nonumber
\end{align}
Thus, by \eqref{48} and \eqref{49}, we have
\begin{align}
\beta_{n,\lambda}^{(k)}&=\sum_{j=0}^{n}\frac{(-1)^{n-j}j!}{(j+1)^{k}}S_{2,\lambda}(n,j),\label{50} \\
\beta_{n,\lambda}^{(k)}(x)&=\sum_{l=0}^{n}\binom{n}{l}(-1)^{n-l}\beta_{l,\lambda}^{(k)}(x)_{n-l,\lambda}.\label{51}
\end{align}
From \eqref{43}, we note that
\begin{align}
&\sum_{n=0}^{\infty}\frac{1}{(n+1)^{r}}\bigg\{\sum_{k=0}^{n}\binom{n}{k}(-1)^{k}f_{\lambda}(y+zk)\bigg\}\label{52} \\
&=\sum_{m=0}^{\infty}a_{m}\sum_{p=0}^{m}\binom{m}{p}z^{p}(y)_{m-p,\lambda}\sum_{j=0}^{p}(-1)^{j}j!S_{2,\frac{\lambda}{z}}(p,j)\frac{1}{(j+1)^{r}}\nonumber \\
&=\sum_{m=0}^{\infty}a_{m}\sum_{p=0}^{m} \binom{m}{p}z^{p}(y)_{m-p,\lambda}(-1)^{p}\beta_{p,\frac{\lambda}{z}}^{(r)}\nonumber \\
&=\sum_{m=0}^{\infty}a_{m}\sum_{p=0}^{m}\binom{m}{p}z^{m}z^{p-m}\bigg(\frac{y}{z}\bigg)_{m-p,\frac{\lambda}{z}}(-1)^{p}\beta_{p,\frac{\lambda}{z}}^{(r)}\nonumber \\
&=\sum_{m=0}^{\infty}a_{m}(-1)^{m}z^{m}\sum_{p=0}^{m}\binom{m}{p}\bigg(\frac{y}{z}\bigg)_{m-p,\frac{\lambda}{z}}(-1)^{m-p}\beta_{p,\frac{\lambda}{z}}^{(r)}\nonumber \\
&=\sum_{m=0}^{\infty}a_{m}(-1)^{m}z^{m}\beta_{m,\frac{\lambda}{z}}^{(r)}\bigg(\frac{y}{z}\bigg).\nonumber
\end{align}
Therefore, by \eqref{52}, we obtain the following theorem.
\begin{theorem}
Let 	$\displaystyle f_{\lambda}(t)=\sum_{k=0}^{\infty}a_{k}(t)_{k,\lambda}\in\mathbb{C}[\![t]\!]\displaystyle$. Then the following identity is valid.
\begin{displaymath}
\sum_{n=0}^{\infty}\frac{1}{(n+1)^{r}}\bigg\{\sum_{k=0}^{n}\binom{n}{k}(-1)^{k}f_{\lambda}(y+zk)\bigg\}= \sum_{m=0}^{\infty}a_{m}(-1)^{m}z^{m}\beta_{m,\frac{\lambda}{z}}^{(r)}\bigg(\frac{y}{z}\bigg),\quad (r\in\mathbb{Z}).
\end{displaymath}
In particular, for $y=0$, we have
\begin{displaymath}
\sum_{n=0}^{\infty}\frac{1}{(n+1)^{r}}\bigg\{\sum_{k=0}^{n}\binom{n}{k}(-1)^{k}f_{\lambda}(zk)\bigg\}= \sum_{m=0}^{\infty}a_{m}(-1)^{m}z^{m}\beta_{m,\frac{\lambda}{z}}^{(r)},\quad (r\in\mathbb{Z}).
\end{displaymath}
\end{theorem}
Let $f_{\lambda}(t)=(t)_{m,\lambda}$, and let $z=1$ in Theorem 13. Then we have
\begin{equation}
\sum_{n=0}^{m}\frac{1}{(n+1)^{r}}\bigg\{\sum_{k=0}^{n}\binom{n}{k}(-1)^{k}(y+k)_{m,\lambda}\bigg\}=(-1)^{m}\beta_{m,\lambda}^{(r)}(y),\label{53}
\end{equation}
in view of Remark 2.  \par
Now, we observe that
\begin{align}
\sum_{n=0}^{\infty}\frac{(-\lambda)^{n}(1)_{n+1,\frac{1}{\lambda}}}{(n+1)!}\sum_{k=0}^{n}\binom{n}{k}(-1)^{k}e_{\lambda}^{k}(z)&=\sum_{n=0}^{\infty}\frac{(-\lambda)^{n}(1)_{n+1,\frac{1}{\lambda}}}{(n+1)!}(1-e_{\lambda}(z))^{n}\label{54}\\
&=\frac{1}{1-e_{\lambda}(z)}\sum_{n=1}^{\infty}\big(1-e_{\lambda}(z)\big)^{n}\frac{(-\lambda)^{n-1}}{n!}(1)_{n,\frac{1}{\lambda}}\nonumber \\
&=\frac{1}{1-e_{\lambda}(z)}\Big(-\log_{\lambda}(1-(1-e_{\lambda}(z)))\Big)\nonumber \\
&=\frac{z}{e_{\lambda}(z)-1}=\sum_{n=0}^{\infty}\beta_{n,\lambda}\frac{z^{n}}{n!}. \nonumber
\end{align}
Thus, we have
\begin{align}
\sum_{n=0}^{\infty}\bigg(\sum_{m=0}^{\infty}\frac{(-\lambda)^{m}(1)_{m+1,\frac{1}{\lambda}}}{(m+1)!}\sum_{k=0}^{m}\binom{m}{k}(-1)^{k}(k)_{n,\lambda}\bigg)\frac{z^{n}}{n!}=\sum_{n=0}^{\infty}\beta_{n,\lambda}\frac{z^{n}}{n!}. \label{55}	
\end{align}
Therefore, by comparing the coefficients on both sides of \eqref{55}, we obtain the following theorem.
\begin{theorem}
For $n\ge 0$, we have the identity.
\begin{align*}
\beta_{n,\lambda}&=\sum_{k=0}^{\infty}\bigg(\sum_{m=k}^{\infty}\binom{m}{k}\frac{(-\lambda)^{m}(1)_{m+1,\frac{1}{\lambda}}}{(m+1)!}(-1)^{k}\bigg)(k)_{n,\lambda}\\
&=\sum_{m=0}^{\infty}\sum_{k=0}^{m}\binom{m}{k}(-1)^{k}\frac{(-\lambda)^{m}(1)_{m+1,\frac{1}{\lambda}}}{(m+1)!}(k)_{n,\lambda}.
\end{align*}
\end{theorem}
From \eqref{8}, we note that
\begin{align}
\sum_{n=0}^{\infty}F_{n,\lambda}\bigg(-\frac{1}{2}\bigg)\frac{t^{n}}{n!}&=\frac{1}{1+\frac{1}{2}(e_{\lambda}(t)-1)}=\frac{2}{e_{\lambda}(t)+1}\label{56}\\
&=2\bigg(\frac{1}{e_{\lambda}(t)-1}-\frac{2}{e_{\lambda}^{2}(t)-1}\bigg)\nonumber \\
&=\frac{2}{t}\bigg(\frac{t}{e_{\lambda}(t)-1}-\frac{2t}{e_{\frac{\lambda}{2}}(2t)-1}\bigg)\nonumber \\
&=\frac{2}{t} \sum_{n=1}^{\infty}\big(\beta_{n,\lambda}-2^{n}\beta_{n,\frac{\lambda}{2}}\big)\frac{t^{n}}{n!}\nonumber \\
&=\sum_{n=0}^{\infty}\frac{2}{n+1}\big(\beta_{n+1,\lambda}-2^{n+1}\beta_{n+1,\frac{\lambda}{2}}\big)\frac{t^{n}}{n!}.\nonumber
\end{align}
In addition, we also have
\begin{align}
\frac{1}{1+\frac{1}{2}(e_{\lambda}(t)-1)}&=\sum_{p=0}^{\infty}	\frac{(-1)^{p}}{2^{p}}(e_{\lambda}(t)-1)^{p}\label{57}\\
&=\sum_{p=0}^{\infty}\frac{(-1)^{p}p!}{2^{p}}\frac{1}{p!}(e_{\lambda}(t)-1)^{p}\nonumber\\
&=\sum_{n=0}^{\infty}\bigg(\sum_{p=0}^{n}\frac{(-1)^{p}p!}{2^{p}}S_{2,\lambda}(n,p)\bigg)\frac{t^{n}}{n!}.\nonumber
\end{align}
Therefore, by \eqref{56} and \eqref{57}, we obtain the following theorem.
\begin{theorem}
For $n\ge 0$, the following identity holds.
\begin{align*}
F_{n,\lambda}\bigg(-\frac{1}{2}\bigg)&=\frac{2}{n+1}\big(\beta_{n+1,\lambda}-2^{n+1}\beta_{n+1,\frac{\lambda}{2}}\big)\\
&=\sum_{p=0}^{n}\frac{(-1)^{p}p!}{2^{p}}S_{2,\lambda}(n,p).
\end{align*}	
\end{theorem}
Let us take $x=\frac{1}{2}$ in Theorem 11. Then we have
\begin{align}
\sum_{n=0}^{\infty}\bigg(\frac{1}{2}\bigg)^{n}\bigg\{\sum_{k=0}^{n}\binom{n}{k}(-1)^{k}f_{\lambda}(zk)\bigg\}&=\sum_{m=0}^{\infty}a_{m}z^{m}F_{m,\frac{\lambda}{z}}\bigg(-\frac{1}{2}\bigg) \label{58}\\
&=\sum_{m=0}^{\infty}a_{m}z^{m}\frac{2}{m+1}\Big(\beta_{m+1,\frac{\lambda}{z}}-2^{m+1}\beta_{m+1,\frac{\lambda}{2^{z}}}\Big)\nonumber.
\end{align}
By \eqref{8}, we get
\begin{align}
\sum_{m=0}^{\infty}F_{m,\lambda}\bigg(\frac{x}{1-x}\bigg)\frac{t^{m}}{m!}&=\frac{1}{1-\frac{x}{1-x}(e_{\lambda}(t)-1)}=\frac{1-x}{1-x-x(e_{\lambda}(t)-1)}\label{59}\\
&=\frac{1-x}{1-xe_{\lambda}(t)}.\nonumber
\end{align}
On the other hand,
\begin{align}
\frac{1-x}{1-xe_{\lambda}(t)}&=(1-x)\sum_{n=0}	^{\infty}x^{n}e_{\lambda}^{n}(t)=(1-x)\sum_{n=0}^{\infty}x^{n}\sum_{m=0}^{\infty}(n)_{m,\lambda}\frac{t^{m}}{m!}\label{60} \\
&=\sum_{m=0}^{\infty}\bigg((1-x)\sum_{n=0}^{\infty}x^{n}(n)_{m,\lambda}\bigg)\frac{t^{m}}{m!}.\nonumber
\end{align}
Therefore, by \eqref{63} and \eqref{64}, we obtain the following theorem.
\begin{theorem}
For $m\ge 0$, the following holds true.
\begin{displaymath}
\frac{1}{1-x}F_{m,\lambda}\bigg(\frac{x}{1-x}\bigg)=\sum_{n=0}^{\infty}x^{n}(n)_{m,\lambda}.
\end{displaymath}	
Let
\begin{equation}
h_{\lambda}(t)=\frac{1}{\mu e_{\lambda}(\gamma t)+1}=\frac{1}{1-(-\mu e_{\lambda}(\gamma t))}=\sum_{n=0}^{\infty}(-\mu)^{n}e_{\lambda}^{n}(\gamma t).\label{61}
\end{equation}
\end{theorem}
Invoking Theorem 16, we obtain
\begin{align}
\bigg(\frac{d}{dt}\bigg)^{m}\bigg(\frac{1}{\mu e_{\lambda}(\gamma t)+1}\bigg)&=\bigg(\frac{d}{dt}\bigg)^{m}\sum_{n=0}^{\infty}(-\mu)^{n}e_{\lambda}^{n}(\gamma t)\label{62} \\
&=\sum_{n=0}^{\infty}(-\mu)^{n}(n)_{m,\lambda}\gamma^{m}e_{\lambda}^{n-m \lambda}(\gamma t) \nonumber \\
&=\gamma^{m}\sum_{n=0}^{\infty}\Big(-\mu e_{\lambda}(\gamma t)\Big)^{n}(n)_{m,\lambda}e_{\lambda}^{-m\lambda}(\gamma t)\nonumber \\
&=\frac{\gamma^{m}}{(1+\lambda \gamma t)^{m}}\sum_{n=0}^{\infty}\big(-\mu e_{\lambda}(\gamma t)\big)^{n}(n)_{m,\lambda}\nonumber \\
&=\frac{\gamma^{m}}{(1+\lambda \gamma t)^{m}}\frac{1}{1+\mu e_{\lambda}(\gamma t)}F_{m,\lambda}\bigg(\frac{-\mu e_{\lambda}(\gamma t)}{1+\mu e_{\lambda}(\gamma t)}\bigg).\nonumber
\end{align}
From \eqref{61} and \eqref{62}, we get the following equation:
\begin{equation}
\bigg(\frac{d}{dt}\bigg)^{m}h_{\lambda}(t)\bigg|_{t=0}=h_{\lambda}^{(m)}(0)=\gamma^{m}\frac{1}{1+\mu}F_{m,\lambda}\bigg(\frac{-\mu}{1+\mu}\bigg),\quad (m\ge 0).\label{63}	
\end{equation}
Thus we have
\begin{align}
h_{\lambda}(t)&=\frac{1}{\mu e_{\lambda}(\gamma t)+1}=\sum_{m=0}^{\infty}\frac{h_{\lambda}^{(m)}(0)}{m!}t^{m}=\frac{1}{1+\mu}\sum_{m=0}^{\infty}\gamma^{m}F_{m,\lambda}\bigg(\frac{-\mu}{1+\mu}\bigg)\frac{t^{m}}{m!}.\label{64}
\end{align} \par
Let us take $\mu=1$ and $\gamma=1$ in \eqref{64}. Then we have
\begin{equation}
\frac{1}{e_{\lambda}(t)+1}=\frac{1}{2}\sum_{m=0}^{\infty}F_{m,\lambda}\bigg(-\frac{1}{2}\bigg)\frac{t^{m}}{m!}. \label{65}	
\end{equation}
As is well known, Carlitz's degenerate Euler polynomials are defined by
\begin{equation}
\frac{2}{e_{\lambda}(t)+1}e_{\lambda}^{x}(t)=\sum_{n=0}^{\infty}\mathcal{E}_{n,\lambda}(x)\frac{t^{n}}{n!}.\label{66}	
\end{equation}
When $x=0$, $\mathcal{E}_{n,\lambda}= \mathcal{E}_{n,\lambda} (0)$ are called the degenerate Euler numbers. \par
From \eqref{65}, we note that
\begin{align}
\frac{2}{e_{\lambda}(t)+1}e_{\lambda}^{x}(t)&=\bigg(\sum_{m=0}^{\infty}F_{m,\lambda}\bigg(-\frac{1}{2}\bigg)\frac{t^{m}}{m!}	\bigg)e_{\lambda}^{x}(t)\label{67}\\
&=\sum_{n=0}^{\infty}\bigg(\sum_{m=0}^{n}\binom{n}{m}F_{m,\lambda}\bigg(-\frac{1}{2}\bigg)(x)_{n-m,\lambda}\bigg)\frac{t^{n}}{n!}.\nonumber
\end{align}
Therefore, by comparing the coefficients on both sides of \eqref{66} and \eqref{67}, we obtain the folloiwnwg theorem.
\begin{theorem}
For $n\ge 0$, the following identity is valid.
\begin{displaymath}
\mathcal{E}_{n,\lambda}(x)=\sum_{m=0}^{n}\binom{n}{m}F_{m,\lambda}\bigg(-\frac{1}{2}\bigg)(x)_{n-m,\lambda}.
\end{displaymath}	
In particular, for $x=0$, we have
\begin{displaymath}
\mathcal{E}_{n,\lambda}=F_{n,\lambda}\bigg(-\frac{1}{2}\bigg)=\frac{2}{n+1}\Big(\beta_{n+1,\lambda}-2^{n+1}\beta_{n+1,\frac{\lambda}{2}}\Big).
\end{displaymath}
\end{theorem}

\section{Further Remarks}

Here we obtain an expression for $\mathcal{E}_{m,\lambda}\bigg(\frac{1}{2}\bigg)$ and a general operational formula \eqref{77}.

Taking $x=\frac{1}{2}$ and $z=1$ in Theorem 11, we have
\begin{equation}
\sum_{n=0}^{\infty}\bigg(\frac{1}{2}\bigg)^{n}\bigg\{\sum_{k=0}^{n}\binom{n}{k}(-1)^{k}f_{\lambda}(y+k)\bigg\}=\sum_{m=0}^{\infty}a_{m}	\sum_{p=0}^{m}\binom{m}{p}(y)_{m-p,\lambda}F_{p,\lambda}\bigg(-\frac{1}{2}\bigg).\label{68}
\end{equation}
By Theorem 17, we get
\begin{equation}
\sum_{n=0}^{\infty}\bigg(\frac{1}{2}\bigg)^{n}\bigg\{\sum_{k=0}^{n}\binom{n}{k}(-1)^{k}f_{\lambda}(y+k)\bigg\}=\sum_{m=0}^{\infty}a_{m}\mathcal{E}_{m,\lambda}(y).\label{69}
\end{equation}
In particular, for $y=\frac{1}{2}$,
\begin{equation}
\sum_{n=0}^{\infty}\bigg(\frac{1}{2}\bigg)^{n}\bigg\{\sum_{k=0}^{n}\binom{n}{k}(-1)^{k}f_{\lambda}\bigg(\frac{1}{2}+k\bigg)\bigg\}=\sum_{m=0}^{\infty}a_{m}\mathcal{E}_{m,\lambda}\bigg(\frac{1}{2}\bigg).\label{70}
\end{equation} \par
Let us take $f_{\lambda}(t)=(t)_{m,\lambda}$. Then, $a_{k}=\delta_{k,m},~~ (k \ge 0)$.
From \eqref{70}, we have
\begin{align}
\mathcal{E}_{m,\lambda}\bigg(\frac{1}{2}\bigg)&=\sum_{n=0}^{\infty}\bigg(\frac{1}{2}\bigg)^{n}\bigg\{\sum_{k=0}^{n}\binom{n}{k}(-1)^{k}\bigg(\frac{1}{2}+k\bigg)_{m,\lambda}\bigg\} \label{71} \\
&=\sum_{n=0}^{m}\bigg(\frac{1}{2}\bigg)^{m}\sum_{k=0}^{n}\binom{n}{k}(-1)^{k}\bigg(\frac{1}{2}+k\bigg)_{m,\lambda},\nonumber	
\end{align}
by invoking Remark 2. \par
Let $D=\frac{d}{dx}$, and let $\displaystyle f(x)=\lim_{\lambda\rightarrow 0}f_{\lambda}(x)=\sum_{n=0}^{\infty}a_{n}x^{n}\displaystyle$. Then we note that
\begin{align}
f(x^{1-\lambda}D)e^{x}&=\sum_{n=0}^{\infty}a_{n}(x^{1-\lambda}D)^{n}\sum_{l=0}^{\infty}x^{l}\frac{1}{l!}\label{72}\\
&=\sum_{n=0}^{\infty}a_{n}\sum_{l=0}^{\infty}\frac{(l)_{n,\lambda}}{l!}x^{l-n\lambda}\nonumber \\
&=\sum_{n=0}^{\infty}a_{n}\bigg(\sum_{l=0}^{\infty}\frac{(l)_{n,\lambda}}{l!}x^{l}e^{-x}\bigg)e^{x}x^{-n\lambda}\nonumber \\
&=\bigg(\sum_{n=0}^{\infty}a_{n}\phi_{n,\lambda}(x)x^{-n\lambda}\bigg)e^{x}.\nonumber	
\end{align}
On the other hand,
\begin{align}
f(x^{1-\lambda}D)e^{x}&=\sum_{n=0}^{\infty}a_{n}(x^{1-\lambda}D)^{n}\sum_{k=0}^{\infty}\frac{x^{k}}{k!}\label{73}\\
&=\sum_{k=0}^{\infty}\frac{1}{k!}\bigg(\sum_{n=0}^{\infty}a_{n}(k)_{n,\lambda}x^{-n\lambda}\bigg)x^{k}.\nonumber
\end{align}
Thus, we note that
\begin{displaymath}
\sum_{k=0}^{\infty}\frac{1}{k!}\bigg(\sum_{n=0}^{\infty}a_{n}(k)_{n,\lambda}x^{-n\lambda}\bigg)x^{k}=e^{x}\sum_{n=0}^{\infty}a_{n}\phi_{n,\lambda}(x)x^{-n\lambda}.
\end{displaymath}
Let $\displaystyle f(t)=\sum_{n=0}^{\infty}a_{n}t^{n},\ g(x)=\sum_{k=0}^{\infty}c_{k}x^{k}\displaystyle$. Then we have
\begin{align}
f(x^{1-\lambda}D)g(x)&=\sum_{n=0}^{\infty}a_{n}(x^{1-\lambda}D)^{n}\sum_{k=0}^{\infty}c_{k}x^{k}\label{74} \\
&=\sum_{k=0}^{\infty}c_{k}\bigg(\sum_{n=0}^{\infty}a_{k}(k)_{n,\lambda}x^{-n\lambda}\bigg)x^{k}.\nonumber
\end{align}
By Taylor expansion, we get
\begin{equation}
\sum_{k=0}^{\infty}\frac{g^{(k)}(0)}{k!}x^{k}=g(x)=\sum_{k=0}^{\infty}c_{k}x^{k}.\label{75}
\end{equation}
Thus, we have $c_{k}=\frac{g^{(k)}(0)}{k!},\ (k\ge 0)$. \par
From \eqref{74} and \eqref{75}, we note that
\begin{align}
f(x^{1-\lambda}D)g(x)&=\sum_{k=0}^{\infty}\frac{g^{(k)}(0)}{k!}\bigg(\sum_{n=0}^{\infty}a_{n}(k)_{n,\lambda}x^{-n\lambda}\bigg)x^{k}\label{76} \\
&=\sum_{n=0}^{\infty}\bigg\{a_{n}\sum_{k=0}^{\infty}\frac{g^{(k)}(0)}{k!}(k)_{n,\lambda}x^{k}\bigg\}x^{-n\lambda}.\nonumber	
\end{align}
On the other hand,
\begin{align}
f(x^{1-\lambda}D)g(x)&=\sum_{n=0}^{\infty}\frac{f^{(n)}(0)}{n!}(x^{1-\lambda}D)^{n}g(x)\label{77} \\
&=\sum_{n=0}^{\infty}\bigg(\frac{f^{(n)}(0)}{n!}\sum_{k=0}^{n}S_{2,\lambda}(n,k)x^{k}D^{k}g(x)\bigg)x^{-n\lambda},\nonumber	
\end{align}
where we used
\begin{align}
(x^{1-\lambda}D)^{n}f(x)=x^{-n \lambda}\sum_{k=0}^{n}S_{2,\lambda}(n,k)x^{k}D^{k}f(x). \label{78}
\end{align}
The operational formula \eqref{78} follows by induction $n$ from the following recurrence relation:
\begin{align*}
S_{2,\lambda}(n+1,k)=S_{2,\lambda}(n,k-1)+(k-n \lambda)S_{2, \lambda}(n,k).
\end{align*}

\section{Conclusion}
In this paper, we adopted the ideas of Boyadzhiev on binomial power series and expressed several identities involving degenerate formal power series as those including degenerate Stirling numbers of the second kind, degenerate Bell polynomials, degenerate Fubini polynomials and degenerate poly-Bernoulli polynomials. In addition, we also obtained some related identities which involve the $\lambda$-falling factorials, the degenerate Stirling numbers of both kinds, the degenerate Bernoulli numbers, the degenerate Fubini polynomials and the degenerate Euler polynomials. \par
Here our replacement of power series by degenerate power series is in the same spirit as the recent paper
[8].The Rota's theory on umbral calculus is based on the linear functionals and the differential operators. The Sheffer sequences occupy the central position in the theory and are characterized by the generating functions involving  the usual exponential function. The motivation for [8] started from the question that what if the usual exponential function is replaced by the degenerate exponential functions. It may be said that this question is very natural in view of the regained recent interests in degenerate special numbers and polynomials. As it turns out, it corresponds to replacing the linear functional by the $\lambda$-linear functionals and the differential operator by the $\lambda$-differential operators. In this way, we were led to introduce $\lambda$-umbral calculus and $\lambda$-Sheffer sequences. \par
As one of our future projects, we would like to continue to pursue our searches for $\lambda$-counterparts of some special polynomials, some special numbers, some transcendental functions and so on.


\begin{thebibliography}{9}
\bibitem{1}
Araci, S. \emph{Novel identities involving Genocchi numbers and polynomials arising from applications of umbral calculus,} Appl. Math. Comput. \textbf{233} (2014), 599--607.
\bibitem{2}
Bayad, A.; Hamahata, Y. \emph{Polylogarithms and poly-Bernoulli polynomials,} Kyushu J. Math. \textbf{65} (2011), no. 1, 15--24.
\bibitem{3}
Boyadzhiev, K. N. \emph{Power series with binomial sums and asymptotic expansion,}
Int. J. Math. Anal. \textbf{8} (2014), 1389--1414.
\bibitem{4}
Boyadzhiev, K. N. \emph{Series transformation formulas of Euler type, Hadamard product of functions, and harmonic number identities,} Indian J. Pure Appl. Math. \textbf{42} (2011), 371--387.
\bibitem{5}
Boyadzhiev, K. N.\emph{ A series transformation formula and related polynomials,} Int. J. Math. Math. Sci. \textbf{2005} (2005), no. 23, 3849-3866.
\bibitem{6}
Carlitz, L. \emph{Degenerate Stirling, Bernoulli and Eulerian numbers,} Utilitas Math. \textbf{15}  (1979), 51--88.
\bibitem{7}
Kim, D. S.; Kim, T. \emph{A note on a new type of degenerate Bernoulli numbers,} Russ. J. Math. Phys. \textbf{27} (2020), no. 2, 227--235.
\bibitem{ }
Kim, D. S.; Kim, T. \emph{Degenerate Sheffer sequences and $\lambda$-Sheffer sequences,} J. Math. Anal. Appl. \textbf{493} (2021), no. 1, Paper No. 124521, 21 pp.
\bibitem{8}
Kim, H. K. \emph{Fully degenerate Bell polynomials associated with degenerate Poisson random variables,} Open Math. \textbf{19} (2021), no. 1, 284--296.
\bibitem{9}
Kim, T.; Kim, D. S.; \emph{Some Identities on Truncated Polynomials Associated with Degenerate Bell Polynomials,} Russ. J. Math. Phys. \textbf{28} (2021), no. 3, 342--355.
\bibitem{10}
Kim, T.; Kim, D. S. \emph{Degenerate Laplace transform and degenerate gamma function.} Russ. J. Math. Phys. \textbf{24} (2017), no. 2, 241--248.
\bibitem{11}
Kim, T.; Kim, D. S. \emph{Note on the degenerate gamma function.} Russ. J. Math. Phys. \textbf{27} (2020),  no. 3, 352--358.
\bibitem{12}
Kim, T.; Kim, D. S.; Dolgy, D. V. \emph{On partially degenerate Bell numbers and polynomials,} Proc. Jangjeon Math. Soc. \textbf{20} (2017),  no. 3, 337--345.
\bibitem{13}
Kim, T.; Kim, D. S.; Jang, G.-W. \emph{A note on degenerate Fubini polynomials,} Proc. Jangjeon Math. Soc. \textbf{20} (2017),  no. 4, 521--531.
\bibitem{14}
Kwon, J.; Kim, W. J.; Rim, S.-H. \emph{On the some identities of the type 2 Daehee and Changhee polynomials arising from $p$-adic integrals on $\mathbb{Z}_{p}$,} Proc. Jangjeon Math. Soc. \textbf{22} (2019), no. 3, 487--497.
\bibitem{15}
Park, J.-W.; Kim, B. M.; Kwon, J. \emph{Some identities of the degenerate Bernoulli polynomials of the second kind arising from $\lambda$-Sheffer sequences,} Proc. Jangjeon Math. Soc. \textbf{24} (2021),  no. 3, 323-342.
\bibitem{16}
Roman, S. \emph{The umbral calculus,} Pure and Applied Mathematics 111, Academic Press, Inc. [Harcourt Brace Jovanovich, Publishers], New York, 1984.
\bibitem{17}
 Simsek, Y. \emph{Identities and relations related to combinatorial numbers and polynomials,}  Proc. Jangjeon Math. Soc.  \textbf{20} (2017), no. 1, 127–135.

\end{thebibliography}
\end{document}